\documentclass[reqno,12pt]{amsart}
\usepackage{latexsym,amssymb,amsfonts, amsmath}
\usepackage{a4wide}
\newtheorem{theorem}{Theorem}

\newtheorem{proposition}[theorem]{Proposition}
\theoremstyle{definition}

\newtheorem{remark}[theorem]{Remark}

\newcommand{\eps}{\varepsilon}

\newcommand{\ee}{\mathbb{E}}

\newcommand{\nn}{\mathbb{N}}
\newcommand{\N}{\mathbb{N}}
\newcommand{\rr}{\mathbb{R}}
\newcommand{\R}{\mathbb{R}}
\newcommand{\pp}{\mathbb{P}}

\newcommand{\ds}{\displaystyle}

\newcommand{\CCalpha}{{\mathcal C}^{\alpha}}
\def\CC{\mathcal C}

\def\BB{\mathcal B}

\begin{document}
\bibliographystyle{plain}

\author{Fran\c cois Bolley}

\address{ENS Lyon, UMPA (UMR 5669), 46 all\'ee d'Italie, F-69364 Lyon Cedex 07}
\address{Current address: Institut de math\'ematiques - LSP (UMR C5583), Universit\'e Paul-Sabatier, Route de Narbonne, F-31062 Toulouse cedex 4}
\email{bolley@cict.fr}

\
\title[]
{Quantitative concentration inequalities on sample path space for mean field interaction}

\thanks{}

\begin{abstract}
We consider a system of particles experiencing diffusion and mean field interaction, 
and study its behaviour when the number of particles goes to infi\-nity. We derive non-asymptotic large deviation bounds 
measuring the concentration of the empirical measure of the paths of the particles around its limit. The method is based on a coupling argument, strong integrability estimates on the paths in H\"older norm, and some~gen\nobreak eral concentration result for the empirical measure of
identically distributed independent paths.
\end{abstract}

\maketitle
\tableofcontents

\section*{Introduction}

\medskip

This paper is devoted to the study of the behaviour of some large stochastic particle system. In the models to be considered, the evolution of each particle is governed 
by a random diffusive term, an exterior force field and a mean field interaction with the other 
particles. For such models the limit behaviour has been clearly identified and studied in terms of law of large numbers, central limit theorem and large deviations. Here we shall give new quantitative estimates on the convergence in the setting of large deviations.

\medskip

This follows some works addressing this issue at the level of observables 
or at the level of the whole system at a given time, that we now summarize. For that purpose let $(X_t^i)_{1 \leq i \leq N}$ be 
the position at time $t$ of the $N$ 
particles in the phase space $\rr^d$ and let $\mu_t$ be some probability measure describing the limit
behaviour of the system. At the level of Lipschitz observables, F. Malrieu \cite{Mal01} adapted ideas of concentration 
of measure to obtain bounds~like
\begin{equation}\label{conc-obs}
\sup_{[\varphi]_1\leq 1} \pp \left[ \Bigl | 
\frac{1}{N} \sum_{i=1}^{N} \varphi (X_t^i) - \int_{\rr^d} \varphi\,d\mu_t \Bigr | > \frac{C}{\sqrt{N}} + \eps \right]
\leq 2\, e^{-\lambda N \eps^2} \, , \qquad N \geq 1,
\end{equation}
where $C$ and $\lambda$ are some constants independent of $\eps$ and $N$, and $[ \, \cdot \, ]_1$ is the Lipschitz seminorm defined by
\[ 
[\varphi]_1 := \sup_{x\neq y} \frac{|\varphi(x)-\varphi(y)|}{\vert x - y \vert} \cdot
\]
 In other words, letting $\delta_x$ stand for the Dirac mass at a point $x \in \rr^d$, the empirical measure
\[
\hat{\mu}_t^N := \frac{1}{N} \sum_{i=1}^N \delta_{X_t^i} 
\]
of the system, which generates the observables at time $t$, satisfies the deviation inequality
\[
\sup_{[\varphi]_1 \leq 1} \pp \left[ \Bigl | 
\int_{\rr^d} \varphi\,d\hat{\mu}_t^N - \int_{\rr^d} \varphi\,d\mu_t \Bigr | > \frac{C}{\sqrt{N}} + \eps \right]
\leq 2\, e^{-\lambda N \eps^2}\, , \qquad N \geq 1. 
\]

Now one can measure how this empirical measure $\hat{\mu}_t^N$ is close to its limit $\mu_t$ in a stronger sense, namely, at the very level of the measures. For this, adapting Sanov's large deviation argument, the authors in \cite{BGV05} got quantitative and non-asymptotic bounds on the deviation of $\hat{\mu}_t^N$ 
around $\mu_t$ for some distance which induces a topology stronger that the narrow topology. By comparison with
\eqref{conc-obs}, these bounds can be written as
\begin{equation}\label{conc-loi}
\pp \left[ \sup_{[\varphi]_{1}\leq 1} \Bigl | 
\frac{1}{N} \sum_{i=1}^{N} \varphi (X_t^i) - \int_{\rr^d} \varphi\,d\mu_t \Bigr | > \eps \right]
\leq C(\eps) \,  e^{- \lambda N \eps^2} \, , \qquad N \geq 1. 
\end{equation}

\bigskip

In this work we want to go one step further by considering the {\bf trajectories} of the particles.
A natural object to consider is the empirical measure of the  trajectories $(X_t^i)_{0\leq t \leq T}$ on some given time 
interval $[0,T]$, which is defined as
\[
\hat{\mu}_{[0,T]}^N := \frac{1}{N} \sum_{i=1}^N \delta_{(X_t^i)_{0\leq t \leq T}} 
\]
where $\delta_{(X_t^i)_{0\leq t \leq T}}$ is the Dirac mass on the path $(X_t^i)_{0\leq t \leq T}$.
This is a random probability measure, no longer on the phase space $\rr^d$, but now on the path space, which in our model is 
the space of $\rr^d$-valued continuous functions on $[0,T]$.

 Its limit behaviour is given as follows: the limit $\mu_t$ of the 
empirical measure $\hat{\mu}_t^N$ at time $t$ can be seen as the law of the solution at time $t$ to a stochastic 
differential equation; then 
the law of the whole process on $[0,T]$ so defined will be the limit of $\hat{\mu}_{[0,T]}^N$. We shall give a precise meaning to
this convergence, and some estimates which are the analogue of \eqref{conc-loi} in the path space; in particular we shall see that they imply \eqref{conc-loi} by projection at time $t$.

\medskip

In the next section we state our main results and give an insight of the proofs, which will be given in more detail in the 
following sections. 

\medskip

\section{Statement of the results}

\smallskip

\subsection{Some notation and definitions}\label{subsectionnotation} 

One of the key points in this work is to measure the discrepancy between probability measures:
this will be done by means of {\bf Wasserstein distances}, which have revealed convenient in this 
type of issues and are defined as follows. Let $(X,d)$ be a separable and complete metric space, and $p$ be a real number $\geq 1$; the Wasserstein distance of order $p$ between two Borel probability measures $\mu$ and 
$\nu$ on $X$  is
\[
W_p(\mu, \nu) := \inf_{\pi} \left( \iint_{ X \times X} d(x,y)^p \, d\pi(x,y) \right)^{1/p}
\]
where $\pi$ runs over the set of all joint measures on $X \times X$ with marginals $\mu$ and $\nu$. $W_p$ induces a metric on the set of Borel probability measures on $X$ with finite moment 
$\ds \int_{X} d(x_0, x)^p \, d\mu(x)$ for some (and thus any) $x_0$ in $X$; convergence in this metric is equi\-valent to narrow convergence
(against bounded continuous functions) plus some tightness condition on the moments (see for instance \cite{Bol05}, \cite{Vil03} for further 
details on these distances). 

\smallskip
At some point the space $(X,d)$ will be $\rr^d$ equipped with the Euclidean distance $ \vert \cdot \vert$, and in
this context $W_p$ will be denoted $W_{p, \tau}$. But $(X,d)$ will mainly be the space ${\mathcal C} ([0,T], \rr^d)$, also denoted
${\mathcal C}$ if no confusion is possible, of $\rr^d$-valued continuous functions on $[0,T]$, equipped with the uniform norm
\[
\Vert f \Vert_\infty := \sup_{0 \leq t \leq T}\vert f (t)\vert; 
\]
for this space $W_p$ will be denoted $W_{p,[0, T]}$. The Wasserstein distances considered 
in these two situations are linked in the following way: if, for $0 \leq t \leq T$, $\pi_t$ is the projection from ${\mathcal C}$ 
into $\rr^d$ defined by $\pi_t(f) = f(t)$, then for any Borel probability measures $\mu$ and $\nu$ on ${\mathcal C}$, and any 
$p \geq 1$, the relation
\begin{equation}\label{ineq-projwp}
W_{p, \tau}(\pi_t \sharp \mu, \pi_t \sharp \nu) \leq W_{p, [0,T]} (\mu, \nu), \qquad 0 \leq t \leq T
\end{equation}
holds, where $\pi_t \sharp \mu$ is the image measure of $\mu$ by $\pi_t$. 

\bigskip

The distance between two probability measures $\mu$ and $\nu$ on $X$ can also be expressed in terms of the {\bf relative entropy} of
$\nu$ with respect to $\mu$ (for instance), defined by
\[
H(\nu | \mu) = \int_{X} \frac{d\nu}{d\mu} \, \ln \frac{d\nu}{d\mu} \, d\mu
\]
if $\nu$ is absolutely continuous with respect to $\mu$, and $H(\nu | \mu)  = +\infty$ otherwise. 

Both notions are linked by the family of {\it transportation} or {\it Talagrand inequalities}: given $p \geq 1$ and $\lambda >0$,
we say that a probability measure $\mu$ on $X$ satisfies the inequality $T_p(\lambda)$ if
\[
W_p(\mu, \nu) \leq \sqrt{\frac{2}{\lambda} \, H (\nu | \mu)}
\]
holds true for any measure $\nu$, and $\mu$ satisfies $T_p$ if it satisfies $T_p(\lambda)$ for some $\lambda >0$. By 
Jensen's inequality, the weakest of all is $T_1$, which is also the only one for which a simple charac\-terization is known: a 
measure $\mu$ satisfies $T_1(\lambda)$ for some $\lambda >0$ if and only if it admits a square-exponential moment, in the sense that 
there exist $a >0$ and $x_0$ in $X$~such that $\ds \int_{X} e^{a d(x_0, x)^2} \, d\mu(x)$ be finite. Numerical relations 
between such $a$ and $\lambda$ are given in~\cite{BV05, DGW04}.

\medskip

\subsection{A general concentration inequality for empirical measures}\label{subsectionprelim}
The proof of our main theorem on the particle system is based on some general concentration result for the empirical measure
of ${\mathcal C}$-valued independent and identically distributed random variables. We state this result separately.

For this purpose, given some Borel probability measure $\mu$ on $\CC$ and $N$ independent random variables 
$(X^i)_{1 \leq i \leq N}$  with law $\mu$, we let
\[
\hat{\mu}^{N} := \frac{1}{N} \sum_{i=1}^{N}  \delta_{X^i}
\]
denote their empirical measure. 

Given some  real number $\alpha \in (0,1]$, we let  
$\CCalpha := {\mathcal C}^{\alpha}([0,T], \R^d)$ be the space of  functions in  $ \CC := {\mathcal C}([0,T], \R^d)$ which  moreover 
are  H\"older of order  $\alpha$, equipped with the H\"older norm
\[
\Vert f \Vert_\alpha :=\sup\, ( \Vert f \Vert_\infty , [f ]_{\alpha} )
\]
where 
\[
[ f ]_{\alpha}  :=  \sup_{0 \leq t,s \leq T}\frac {\vert f(t)-f(s)\vert }{\vert t-s \vert^\alpha} \cdot
\]
$\CCalpha$ is a Borel set of the space $\CC$ equipped with the topology induced by the uniform norm, and for Borel measures 
on $\CC$, concentrated on $\CCalpha$, we have in the above notation:

\medskip

\begin{theorem}\label{thmprelim}
Let $p \in [1,2]$ and let $\mu$ be a Borel probability measure on $\CC$ satisfying a 
$T_p(\lambda)$ inequality for some $\lambda >0$, and such that $\ds \int_{\CC} e^{ a \Vert x \Vert^2_{\alpha}} \, d\mu(x)$ be finite 
for some $a >0$ and $\alpha \in (0,1]$. Then, for any $\alpha'< \alpha$ and $\lambda' < \lambda$, there exists some constant $N_0$ such that 
\begin{equation}\label{conc-Tp2} 
\pp\left[ W_{p, [0, T]}(\mu,\hat{\mu}^N)>\eps\right]\le e^{- \beta_p \, \frac{\lambda'}{2} \, N \, \eps^2}
\end{equation}
for any $\eps >0$ and $N \geq N_0 \, \eps^{-2} \, \exp \, (N_0 \, \eps^{- 1/ \alpha'})$, where
$$
\beta_p = \left\{ \begin{array}{clc}
                           1 & \mbox{if} & 1\leq p < 2\\
                            (1 + \sqrt{\lambda / a})^{-2} & \mbox{if} & p=2.   
                          \end{array}
                  \right.
$$
\end{theorem}

\medskip

Here the constant $N_0$ depends on $\mu$ only through $\lambda, a, \alpha$ and 
$\ds \int_{\CC} e^{a \Vert x \Vert^2_{\alpha}} \, d\mu(x)$.

\smallskip

Let us make a few remarks on this result.  

\smallskip

First of all, for another formulation of the obtained bound in the case when $p=1$, we recall Kantorovich-Rubinstein dual expression of the $W_1$ distance 
on a general space $(X, d)$:
\begin{equation}\label{KR}
W_1(\mu, \nu) = \sup_{[\varphi]_1 \leq 1} \Big\{ \int_{X} \varphi \, d\mu - \int_{X} \varphi \, d\nu \Big\}
\end{equation}
where $ \ds [\varphi]_{1} := \sup_{x \neq y} \frac{\vert \varphi(x)-\varphi(y) \vert}{d(x,y)} \cdot$
Then a result by S. Bobkov and F. G\"otze \cite{BG99} ensures that a $T_1(\lambda)$ inequality for $\mu$ is equivalent to 
the concentration inequality
\[
\sup_{[\varphi]_{1}\leq 1} \pp \left[ \frac{1}{N} \sum_{i=1}^{N} \varphi (X^i) - \int_{X} \varphi\,d\mu  > \eps \right]
\leq  e^{-\frac{\lambda}{2} N \eps^2} \, , \qquad N \geq 1 \, .
\]
By comparison, the bound given by Theorem \ref{thmprelim} implies
\[
\pp \left[  \sup_{[\varphi]_{1}\leq 1} \Bigl( 
\frac{1}{N} \sum_{i=1}^{N} \varphi (X^i) - \int_{\CC} \varphi\,d\mu \Bigr) > \eps \right]
\leq  e^{-\frac{\lambda'}{2} N \eps^2} \, , \quad \lambda' < \lambda \, , \quad N \, \text{large enough} \, 
\]
by \eqref{KR}, but a modification of the proof would also lead to
\begin{equation}\label{conc-modif}
\pp \left[  \sup_{[\varphi]_{1}\leq 1} \Bigl( 
\frac{1}{N} \sum_{i=1}^{N} \varphi (X^i) - \int_{\CC} \varphi\,d\mu \Bigr) > \eps \right]
\leq C(\eps) \, e^{-\frac{\lambda'}{2} N \eps^2}, \quad \lambda' < \lambda \, , \quad N \geq 1,
\end{equation}
for some computable large constant $C(\eps)$. In other words we control a much stronger quantity, up to some
loss on the constant in the right-hand side, or some condition on the size of the sample.

\medskip

This result seems reasonable in view of Sanov's theorem (stated in \cite{DZ98} for instance). By applying this theorem to $A := \{ \nu; W_{p,  [0, T]} (\nu, \mu) \geq \eps \}$, for some
given $\eps >0$, one can hope for an upper bound like
\[
\pp \, \big[ W_{p,  [0, T]}(\mu, \hat{\mu}^N) \geq \eps \big] \ds \leq \exp \Big( - N \inf \big\{ H(\nu | \mu) ; \nu \in A \big\} \Big)
\]
for large $N$. With this bound in hand, since
\[
 \inf \big\{ H(\nu | \mu) ; \nu \in A \big\} \geq \frac{\lambda}{2}  \, \eps^2
\]
as $\mu$ satisfies $T_p(\lambda)$, one indeed obtains an upper bound like \eqref{conc-Tp2}, but only in an asymptotic way, whereas Theorem \ref{thmprelim} moreover gives an estimate on a sufficient size of the sample for the deviation bound to hold. Sanov's theorem does not actually
give such an upper bound here; indeed, on an unbounded space such as $\CC$, the closure 
$\overline{A}$ of $A$ (for the narrow topology) contains $\mu$ itself: in particular 
$\inf  \{ H(\nu | \mu) ; \, \nu \in \overline{A} \} = 0 $ and Sanov's theorem only gives the trivial upper bound $\pp \, [ A ] \leq \exp \big( - N \inf \{ H(\nu | \mu) ; \nu \in \overline{A} \} \big) = 1$. 

\smallskip

To get a more relevant upper bound we impose some extra integrability assumption that may at first sight seem strong and odd. 
The reason is that, proceeding as in \cite{BGV05}, we first reduce the issue to a compact set of $\CC$
that has almost full $\mu$ measure: a large ball of $\CCalpha$ will do by Ascoli's theorem and the integrability assumption on $\mu$.
Then on this compact set one can get precise upper bounds by using some techniques based on a covering argument and developed
in \cite{KZ95} (see also \cite[Exercises 4.4.5 and 6.2.19]{DZ98} and \cite{Goz05}).
And actually the assumption is satisfied by the Wiener measure on $\CC$ (recall that the Brownian motion paths are almost
surely H\"older of order $\alpha$ for any $\alpha < 1/2$) and by extension by the
law of the process to be considered.

This integrability assumption again implies the existence of a square-exponential moment for $\mu$ on $\CC$ (for
the uniform norm). Since this is equivalent to some $T_1$ inequality for $\mu$, the $T_p(\lambda)$ assumption is redondant 
 when $p=1$ if one does not care of the involved constants.

Finally this result can be seen as an extension of the following similar concentration result given in \cite[Theorem 1.1]{BGV05}
in the case of measures on $\rr^d$: if $m$ satisfies $T_p(\lambda)$, then
\begin{equation}\label{rd}
\pp \, [ W_{p, \tau}(m, \hat{m}^{N}) > \eps ] \leq e^{ - \gamma_p \frac{\lambda'}{2} \, N \, \eps^2}, 
\quad \eps>0 \, , \quad N \geq N_0 \max(\eps^{-(d'+2)}, 1).
\end{equation}
Let indeed $m$ be such a measure on $\rr^d$. Then the law $\mu$ of a constant process on $[0,T]$ initially distributed according to 
$m$ satisfies the assumptions of Theorem \ref{thmprelim} (one can take any $a < \lambda / 2$), and the bound \eqref{rd}
follows by \eqref{ineq-projwp} with the constant $\gamma_p$ obtained in \cite{BGV05}. Note however that the required size of the 
sample is here much larger for small $\eps$.

\medskip

Theorem \ref{thmprelim} will be proved in Section \ref{sectionprelim}.

\smallskip

\subsection{Interacting particle systems}

We now turn to the study of a system of $N$ stochastic interacting particles which positions $X^i_t$ in the phase space $\rr^d$
$(1 \leq i \leq N)$ evolve according to the system of coupled stochastic differential equations
\begin{equation}\label{edsXi}
dX_t^i = \sqrt{2} \, dB_t^i - \nabla V(X_t^i) dt - \frac{1}{N} \sum_{j=1}^{N}  \nabla W(X_t^i-X_t^j) dt,
\qquad 1 \leq i \leq N.
\end{equation}
Here the $B^i$'s are $N$ standard independent Brownian motions on $\rr^d$, $V$ and $W$ are exterior and interaction potentials.

The state of the system at some given time $t$ is given by the observables $ \ds \frac{1}{N} \sum_{i=1}^{N} \varphi(X^i_t)$, and thus can be described by the random probability measure 
\[
\hat{\mu}^N_t := \frac{1}{N} \sum_{i=1}^{N} \delta_{X^i_t}
\]
on the phase space $\rr^d$, called the empirical measure of the system. Under some regularity and growth assumptions on the potentials $V$ and $W$,
and if the particles are initially distributed in a chaotic way, for instance as independent
and identically distributed variables, then $\hat{\mu}^N_t$ converges  to a solution at time $t$ to the
partial differential equation
\begin{equation}\label{edp}
\frac{\partial \mu_t}{\partial t} = \Delta \mu_t +  \nabla \cdot \bigl( \mu_t \nabla(V+W*\mu_t) \bigr)
\end{equation}
as the number $N$ of particles goes to infinity.
This nonlinear diffusive equation, in which $\Delta, \nabla \cdot$ and $\nabla$ respectively denote the Laplace, divergence and gradient
operators in $\rr^d$, is a McKean-Vlasov equation and has been used in \cite{BCCP98} in the modelling of one-dimensional
gra\-nular media.  The convergence of $\hat{\mu}^N_t$ is strongly linked with 
the phenomenon of propagation of chaos for the interacting particles, and both issues have been studied by H.~Tanaka \cite{Tan84}, 
A.-S.~Sznitman \cite{Szn91}, S.~M\'el\'eard~\cite{Mel96} or S.~Benachour, B.~Roynette, D.~Talay and P.~Vallois \cite{BRTV98, BRV98} for instance. Then quantitative estimates on this convergence have been~obtained by F.~Malrieu 
\cite{Mal01} at the level of observables, and later at the very level of the law in \cite{BGV05}. 

\bigskip

In this work we go one step further and study the limit behaviour of the empirical~measure
\[
\hat{\mu}^N_{[0,T]} = \frac{1}{N} \sum_{i=1}^N \delta_{X^i}.
\]
of the trajectories $X^i = (X^i_t)_{0 \leq t \leq T}$ of the particles on some time interval $[0,T]$.

For that purpose, let $Y = (Y_t)_{0 \leq t \leq T}$ be a solution to the stochastic differential equation
\begin{equation}\label{edsY}
dY_t = \sqrt{2} \, dB_t -\nabla V(Y_t) \, dt - \nabla W * \nu_t(Y_t) \, dt
\end{equation}
starting at $Y_0$ distributed according to the initial condition $\mu_0$ in \eqref{edp}, where $\nu_t$ is the law of $Y_t$ at time $t$.
Then, by It\^o's formula, $\nu_t$ also is a solution to equation \eqref{edp} with initial datum $\mu_0$, and a uniqueness result ensures that
actually $\nu_t = \mu_t$. In other words the limit behaviour $\mu_t$ of $\hat{\mu}^N_t$ is the time-marginal
of the law $\mu_{[0,T]}$ of the process $Y$, in the sense that it is the image measure
of $\mu_{[0,T]}$ by the canonical projection $\pi_t$ defined on the path space $\CC$ by $\pi_t(f) = f(t)$.
Since $\hat{\mu}^N_t$ is the time-marginal of the considered empirical measure $\hat{\mu}^N_{[0,T]}$, it is
natural to hope that $\hat{\mu}^N_{[0,T]}$ converge (in some sense) to $\mu_{[0,T]}$. This convergence has indeed been
proved in the first works mentionned above, and here we want to extend to this new setting the techniques developed in \cite{Mal01} and more particularly in \cite{BGV05}.

\medskip

We shall assume that the potentials $V$ and $W$ are twice differentiable on $\rr^d$, with bounded 
hessian matrices in the sense that there exist real numbers $\beta, \beta', \gamma$ and $\gamma'$ such~that
\begin{equation}\label{D2V,D2W}
\beta \, I \leq D^2V (x) \leq \beta' \, I, \qquad \gamma \, I \leq D^2W(x) \leq \gamma' \, I, \qquad  x \in \rr^d.
\end{equation}
In other words the force fields $\nabla V$ and $\nabla W$ are assumed to be Lipschitz on the
whole $\rr^d$. 

\smallskip

Under these assumptions, global existence and uniqueness, pathwise and in law, of the solutions to \eqref{edsXi} and \eqref{edsY} are proven in \cite{Mel96} for instance for square integrable initial
data; moreover the paths are continuous (in time). We shall also assume that the potential $W$, which gives rise to an interaction term,
is symmetric in the sense that $W(-z) = W(z)$ for all $z \in \rr^d$. Then we shall prove

\smallskip

\begin{theorem}\label{thmprincipal}
Let $\mu_0$ be a probability measure on $\rr^d$, admitting a finite square-exponential
moment in the sense that there exists $a_0 >0$ such that
$\ds \int_{\rr^d}  e^{a_0 |x|^2}\,d\mu_0(x)$ be finite.
Let $(X_0^i)_{1 \leq i \leq N}$ be $N$ independent random variables with common 
law $\mu_0$. Given $T \geq 0$, let $(X^i)_{i}$ be the solution of~\eqref{edsXi} on $[0,T]$ with initial value
$(X_0^i)_i$, where $V$ and $W$ are assumed to satisfy~\eqref{D2V,D2W}; let also $\hat{\mu}^N_{[0,T]}$ be 
the empirical measure associated with the $N$ paths $X^i$.
Let finally $\mu_{[0,T]}$ be the law of the process solution of~\eqref{edsY} for some initial value 
distributed according to $\mu_0$. 

Then, for any $\alpha \in (0, 1/2)$, there exist some positive constants $K$ and $N_0$ such that 
\[
\pp \, \left[ W_{1, [0, T]} (\mu_{[0,T]} , \hat{\mu}_{[0,T]}^N) > \varepsilon\right] \leq 
e^{- K \, N\, \eps^2}
\]
for all $\eps>0$ and $N \geq N_0\, \eps^{-2} \, \exp \, (N_0 \, \eps^{-1 / \alpha})$.
\end{theorem}

Here the constants $K$ and $N_0$ depend on $T, V, W, \alpha $ and $\ds \int_{\rr^d}  e^{a_0 |x|^2}\,d\mu_0(x)$.

By Kantorovich-Rubinstein formulation again, this bound can be written as
\begin{equation}\label{conctraj}
\pp \, \Big[ \sup_{[ \varphi ]_1 \leq 1} \Big( \frac{1}{N} \sum_{i=1}^{N} \varphi(X^i) - 
\int_{\CC} \varphi(x) \, d\mu_{[0,T]}(x) \Big) > \eps \Big] \leq e^{- K \, N\, \eps^2} .
\end{equation}
By projection at time $t$, it implies
concentration inequalities for the time-marginals of the empirical measures similar to inequalities \eqref{conc-obs} and even
\eqref{conc-loi}. But above all it gives concentration estimates at the level of the paths. In return we impose some 
stronger condition on the required size of the sample (however, by \eqref{conc-modif}, one can also get less precise estimates
valid for any number $N$ of particles).

Assume for instance that one is interested in the behaviour of a point $Y_t$ evolving according to \eqref{edsY}. Then, from \eqref{conctraj}, one can derive error bounds in the approximation by $\ds \frac{1}{N} \sum_{i=1}^{N} \varphi (X^i)$ of the expectation of quantities $\varphi (Y)$ which depend on the whole path, such as the distance $d (Y, A) = \inf \big\{ \vert Y_t - y \vert; t \in [0, T], y \in A \big\}$ of the trajectory to a given set $A$ in $\rr^d$, which measures how close $Y_t$ has been to $A$, or the maximal distance $\sup \big\{ \vert Y_t - x \vert; t \in [0, T] \big\}$ to a given point $x$ in the phase space $\rr^d$: for instance, under the assumptions of Theorem \ref{thmprincipal}, for any $T \geq 0$ and $\alpha \in (0, 1/2)$ there exist some positive constants $K$ and $N_0$ such that 
\[
\pp \, \Big[ \Big\vert \ee \, [ d (Y, A) ] - \frac{1}{N} \sum_{i=1}^{N} d (X^i, A) \Big\vert > \eps \Big] \leq 
e^{- K \, N\, \eps^2}
\]
for any Borel set $A$ in $\rr^d$, $\eps>0$ and $N \geq N_0\, \eps^{-2} \, \exp \, (N_0 \, \eps^{-1 / \alpha})$.

\medskip

Theorem \ref{thmprincipal} will be proven in detail in Sections \ref{sectioncouplage} and \ref{sectionconclusion} along the following lines.
Following \cite{Szn91} we proceed by coupling, by introducing a family of
$N$ identically distributed processes $Y^i = (Y_t^i)_{0 \leq t \leq T}$ solution to the (nonlinear) stochastic differential equations
$$
\left\{ \begin{array}{rcl}
dY^i_t & = & \sqrt{2} \, dB^i_t - \nabla V(Y^i_t) \, dt - \nabla W*\mu_t (Y^i_t) \, dt \\
Y^i_0 & =  & X^i_0 
	\end{array}
	\right.\qquad 1 \leq i \leq N \, ;
$$
here $\mu_t$ is the solution at time $t$ to \eqref{edp}, but is also the law on $\rr^d$ of any $Y^i_t$ by It\^o's formula, and, for each $i$, $B^i = (B^i_t)_{0 \leq t \leq T}$ is the Brownian motion driving the 
evolution of $X^i$. In particular the paths $Y^i$ are close to the paths $X^i$ and one can prove that 
there exists some constant $C$ (depending only on $T$) such that
\[
W_{1, [0, T]} (\mu_{[0,T]}, \hat{\mu}^N_{[0,T]}) \leq C \, W_{1, [0,T]} (\mu_{[0,T]}, \hat{\nu}^N_{[0,T]})
\]
hold almost surely, where $\ds \hat{\nu}^N_{[0,T]} := \frac{1}{N} \sum_{i=1}^N \delta_{Y^i}$; hence
controlling the distance between  $\mu_{[0,T]}$ and $\hat{\mu}^N_{[0,T]}$ reduces to the same issue with $\mu_{[0,T]}$ and $\hat{\nu}^N_{[0,T]}$. 

But, by definition, the $N$ processes $Y^i$  for $1 \leq i \leq N$ are independent and
distributed according to $\mu_{[0,T]}$. Then Theorem \ref{thmprelim} ensures good concentration estimates for the empirical measure
$\ds \hat{\nu}^N_{[0,T]}$ around the common law $\mu_{[0,T]}$. 
In the end we obtain the bound
\[
\pp \, \big[ W_{1,  [0, T]} (\mu_{[0,T]} , \hat{\mu}_{[0,T]}^N) > \eps \big] \leq 
\pp \, \left[ W_{1,  [0, T]} (\mu_{[0,T]} , \hat{\nu}_{[0,T]}^N) > \frac{\eps}{C} \right] \leq
e^{ - K \, N\, \eps^2 }
\]
under some condition on $\eps$ and $N$. 

The proof is actually an adaptation of the argument given in \cite[Section 1.6]{BGV05} of estimates \eqref{conc-loi} for time-marginals. The current proof turns out to be simpler in the sense that
it consists in  fewer steps; in return each of these steps is somehow more delicate: for instance, as we shall see in the following sections, the proof of Theorem \ref{thmprincipal} requires the computation of the metric entropy of some space of H\"older-continuous functions, and checking that the law of $Y$ fulfills the assumptions of this theorem needs some strong integrability in H\"older norm on solutions to stochastic differential equations.

\medskip

An adaptation of this proof leads to quantitative estimates on the phenomenon of {\bf propa\-gation of chaos}. For instance, letting  
\[
\hat{\mu}^{N,2}_{[0,T]} := \frac{1}{N(N-1)} \sum_{i\neq j} \delta_{(X^i, \, X^j)}
\]
be the empirical measure on  pairs of paths, the asymptotic independence of two paths (among $N$) can be estimated as in

\begin{theorem} 
With the same notation and assumptions as in Theorem \ref{thmprincipal}, for all $T \geq 0$ and $\alpha \in (0, 1/2)$ there exist
some positive constants $K$ and $N_0$ such that 
\[  
\pp \, \left[ W_{1, [0, T]}(\mu_{[0,T]} \otimes \mu_{[0,T]} , \hat{\mu}_{[0,T]}^{N,2}) > \eps \right] 
\leq e^{ - K \, N\, \eps^2 }
 \]
for all $\eps >0$ and $N \geq N_0 \, \eps^{-2} \, \exp( N_0 \, \eps^{-1 / \alpha})$.
\end{theorem}

Here the constants $K$ and $N_0$ depend on $T, V, W, \alpha$ and a finite square-exponential moment of $\mu_0$,
and $W_{1, [0,T]}$ stands for the Wasserstein distance of order $1$ on the pro\-duct space $\CC \times \CC$. 
The proof consists in writing
the coupling argument for pairs of paths and comparing $\mu_{[0,T]} \otimes \mu_{[0,T]}$ and 
$\ds \hat{\nu}^{N,2}_{[0,T]} := \frac{1}{N(N-1)} \sum_{i\neq j} \delta_{(Y^i, Y^j)}$
by means of $\ds \frac{1}{N^2} \sum_{i,  j} \delta_{(Y^i, Y^j)}$.

\bigskip

Let us finally note that it would be desirable to relax the assumptions made on the potentials $V$ and $W$, in particular so as to 
include the interesting case of the cubic potential $W(z) = \vert z \vert^3 / 3$ on $\rr$, which models the interaction among 
one-dimensional granular media (see \cite{BCCP98}). It could also be interesting to consider the whole trajectories $(X_t)_{t \geq 0}$, and derive concentration bounds on functionals such as hitting times for instance.

\bigskip

Before turning to the proofs we briefly recall the {\bf plan of the paper}. In the coming section we prove Theorem \ref{thmprelim}
for general $\CC$-valued independent variables. The study of the particle system is addressed in the following two sections: in
Section 3 we reduce our concentration issue on interacting particles to the same issue for independent variables by a coupling argument,
whereas in Section 4 we check that we can apply our general concentration result to these independent variables; with this in 
hand we can prove Theorem \ref{thmprincipal}. An appendix is devoted to a general metric entropy estimate in a space of H\"older-continuous functions, which
enters the proof of Theorem \ref{thmprelim}.

\medskip

\section{A preliminary result on independent variables}\label{sectionprelim}

The aim of this section is to prove Theorem \ref{thmprelim} for $N$ independent and identically distributed random
variables valued in ${\CC}$. We have seen how this result, applied to the artificial processes $(Y^i_t)_{0 \leq t \leq T}$, enters 
the study of our interacting particle system.

\smallskip

The proof goes in three steps: truncation to a ball $\BB_{R}^{\alpha}$ of $\CCalpha$, compact for the topology induced by
the uniform norm; covering of $\BB_{R}^{\alpha}$ and then of ${\mathcal P}(\BB_{R}^{\alpha})$ by small balls on which one develops
Sanov's argument; conclusion of the proof by optimizing the introduced parameters. Since the argument follows the lines of the proof given in \cite[Section 2.1]{BGV05} in the finite dimensional case, in which $\mu$ is a measure on $\rr^d$, we shall only sketch it, stressing only the bounds specific to our new framework. We refer to \cite{BGV05} for further details.

\medskip

{\bf Step 1. Truncation}.
Given $R > 0$, to be chosen later on, we denote $\BB_{R}^{\alpha}$ the ball  $\{ f \in \CCalpha ; \Vert f \Vert_{\alpha} \leq R\}$ 
of center $0$ and radius $R$ in $\CCalpha$. This set $\BB_{R}^{\alpha}$ is a compact subset of $\CC$ for the topology
induced by the uniform norm $ \Vert \cdot \Vert_{\infty}$: indeed it is relatively compact in $\CC$ by Ascoli's theorem, and closed 
since if $f$ in 
$\CC$ is the uniform limit of a sequence $(f_n)_n$ in $\CCalpha$, then $\ds \Vert f \Vert_{\alpha}  \leq \liminf_{n} 
\Vert f_n \Vert_{\alpha}$, and in particular $f$ belongs to $\BB_{R}^{\alpha}$ if so do the $f_n$.

\smallskip

Letting ${\bf 1}_{\BB_{R}^{\alpha}}$ be the indicator function of $\BB_{R}^{\alpha}$, we truncate $\mu$ into a probability measure $\mu_R$ 
on the ball $\BB_{R}^{\alpha}$, defined as
\[
\mu_R := \frac{{\bf 1}_{\BB_R^{\alpha}} \, \mu}{\mu[\BB_{R}^{\alpha}]} \cdot
\]
Note that $\mu[\BB_{R}^{\alpha}]$ is positive for $R$ larger than some $R_0$ depending only on $a$ and
$\ds E_a := \int_{\CC} e^{a \Vert x \Vert_{\alpha}^2} \, d\mu(x)$. In this step we reduce the concentration problem for $\CC$ to the same issue for the compact ball $\BB_{R}^{\alpha}$, by bounding the quantity 
$\pp \big[ W_p(\mu, \hat{\mu}^N) > \eps \big]$ in terms of $\mu_R$ and an associated empirical measure 
$\ds \hat{\mu}^N_R:=\frac{1}{N}\sum_{k=1}^N \delta_{X^k_R}$ where the $X^k_R$ are independent variables with law $\mu_R$.

\smallskip

Bounding by above the $\Vert \cdot \Vert_{\infty}$ norms by $\Vert \cdot \Vert_{\alpha}$ norms when necessary, we proceed exactly as in  \cite[proof of Theorem 1.1]{BGV05} to obtain the bound 
\begin{multline}\label{concl-trunc-p}
\pp \big[ W_{p,  [0, T]} (\mu,\hat{\mu}^N)>\eps \big] 
\leq \pp \Big[ W_{p,  [0, T]} (\mu_R,\hat{\mu}^N_R)>\eta\, \eps -2 \, E_{a}^{1/p} R  \, e^{-\frac{a}{p}R^2}\Big] \\
 + \exp \left( -  N \big(  \theta (1-\eta)^p \eps^p - E_{a}\, e^{(a_1 - a)\, R^2} \big) \right);
\end{multline}
here $p$ is any real number in  $[1,2)$, $\eta$ in $(0,1)$, $\eps, \theta >0, \,  a_1 < a$ and $R$ is constrained to be larger than
$R_2 \max ( 1 , \theta^{\frac{1}{2-p}} )$
for some constant $R_2$ depending only on $E_a, a, a_1$ and $p$. 

\smallskip

In the case when  $p=2$, we  obtain
\begin{multline}\label{concl-trunc-2}
\pp \big[ W_{2, [0,T]} (\mu,\hat{\mu}^N)>\eps \big] 
\leq \pp\left[W_{2,  [0, T]} (\mu_R,\hat{\mu}^N_R)>\eta\, \eps -2 \, E_{a}^{1/2} R \, e^{-\frac{a}{2}R^2}\right] \\
 + \exp \left( -  N \big(  \frac{a_1}{2} (1-\eta)^2 \eps^2 - 2 E_{a}^2 \, e^{(a_1 - a)\, R^2} \big) \right).
\end{multline}

\bigskip

{\bf Step 2. Sanov's argument on small balls}.
In view of \eqref{concl-trunc-p} for $p < 2$ or \eqref{concl-trunc-2} for $p=2$, we now aim at bounding 
$\pp \, [ \hat{\mu}^{N}_{R} \in {\mathcal A}]$ where
\[
{\mathcal A} := \Bigl\{\nu\in {\mathcal P}(\BB_R^{\alpha});\quad W_{p, [0, T]} (\nu,\mu_R) \geq \eta \, \eps - 
2 \, E_{a}^{1/p} R \, e^{-\frac{a}{p}R^2}\Bigr\}.
\]
For that purpose, reasoning as in \cite{BGV05}, we let $\delta > 0$ and cover ${\mathcal A}$ with ${\mathcal N}({\mathcal A},
\delta)$ balls $(B_i)_{1 \leq i \leq {\mathcal N}({\mathcal A}, \delta)}$ with radius $\delta / 2$ in $W_{p, [0,T]} $ distance. Then one
can develop Sanov's argument on each of these compact and convex balls, and obtain the bound
\begin{equation}\label{cov}
\ds \pp[ \hat{\mu}^{N}_{R} \in {\mathcal A}] \leq \pp \Big[ \hat{\mu}^{N}_{R} \in  \bigcup_{i=1}^{{\mathcal N}({\mathcal A},
\delta)} B_i \Big] \leq \sum_{i=1}^{ {\mathcal N}({\mathcal A}, \delta)} \pp \left[ \hat{\mu}^{N}_{R} \in B_i \right]
\leq \sum_{i=1}^{ {\mathcal N}({\mathcal A}, \delta)} \exp \Big( -  N \inf_{\nu \in B_i} H(\nu \vert \mu_R) \Big).
\end{equation}

Then, from the $T_p(\lambda)$ inequality for $\mu$,  one establishes an approximate $T_p(\lambda)$ inequality for $\mu_R$: namely, for any $\lambda_1 < \lambda$ there exists $K_1$ 
such that 
\[
H(\nu, \mu_R) \geq \frac{\lambda_1}{2} \, W_{p, [0, T]} (\nu, \mu_R)^2 - K_1 \, R^2 \, e^{- a \, R^2}
\]
for any measure $\nu$ on $\BB_R^{\alpha}$. With this inequality in hand, given $1 \leq p <2$ and $\lambda_2 < \lambda_1 < \lambda$,
one deduces from \eqref{cov} the existence of some positive constants $\delta_1, \eta_1$ and $K_1$ such that
\begin{equation}\label{cov2}
\pp\left[W_{p, [0, T]}(\mu_R,\hat{\mu}^N_R)>\eta\, \eps -2 \, E_{a}^{1/p} R e^{-\frac{a}{p}R^2}\right]
\leq {\mathcal N}({\mathcal A}, \delta) \, \exp \left( -  N \Big( \frac{\lambda_2}{2} \, \eps^2 - K_1 R^2 e^{- a R^2} \Big) \right)
\end{equation}
where we have chosen $\delta := \delta_1 \eps$ and $\eta := \eta_1$.

In the case when $p=2$, we do not choose $\eta$ at this stage, and simply obtain
\[
\pp\left[W_{2, [0, T]} (\mu_R,\hat{\mu}^N_R)>\eta\, \eps -2 \, E_{a}^{1/2} R \, e^{-\frac{a}{2}R^2}\right]
\leq {\mathcal N}({\mathcal A}, \delta) \, \exp \left( -  N \Big( \frac{\lambda_2}{2} \, \eta^2 \eps^2 - K_1 R^2 e^{- a R^2} \Big) 
\right)
\]
where $\delta := \delta_1 \eps$.

Then, since ${\mathcal A}$ is a subset of ${\mathcal P}(\BB_{R}^{\alpha})$, Theorem \ref{thm-annexe}  in the Appendix enables to bound ${\mathcal N}({\mathcal A}, \delta)$ with $\delta = \delta_1 \, \eps$ by
\begin{equation}\label{cov3}
\exp \left( K_2 (R \, \eps^{-1})^d \, 3^{K_2 (R \, \eps^{-1})^{1 / \alpha}} \, \ln \left(\max(1, K_2 \, R \, \eps^{-1}) \right) \right)
\end{equation}
for some constant $K_2$ depending neither on $\eps$ nor on $R$.

\smallskip

\begin{remark}
The order of magnitude of this covering number in an infinite-dimensional setting constitutes a main change by comparison with the finite-dimensional setting of \cite{BGV05}, and will influence the final condition on the size $N$ of the sample.
\end{remark}

\bigskip

{\bf Step 3. Conclusion of the argument}.
We first focus on the case when $p \in [1,2)$. Collecting estimates \eqref{concl-trunc-p}, \eqref{cov2} and \eqref{cov3}, 
we obtain, given $\lambda_2 < \lambda$ and $a_1 < a$, the existence of  positive constants $K_1, K_2, K_3$ and $R_3$ depending
on $E_a, a, a_1, \alpha, \lambda$ and $\lambda_2$ such that
\begin{multline}\label{2termsp}
 \pp \left[W_{p, [0, T]} (\mu,\hat{\mu}^N) > \eps \right]  \\
 \leq 
\exp \left( K_2 (R \, \eps^{-1})^d \, 3^{K_2 (R \, \eps^{-1})^{1 / \alpha}} \, \ln \left(\max(1, K_2 \, R \,  \eps^{-1}) \right)
- N \Big( \frac{\lambda_2}{2} \, \eps^2 - K_1 R^2e^{-\alpha R^2}\Big) \right)  \\
   + \exp \left( -  N \big(K_3 \, \theta \, \eps^p - K_4 e^{(a_1 - a)\, R^2} \big) \right) 
\end{multline}
for all $\eps , \theta >0$ and $R \geq R_3 \max (1, \theta^{\frac{1}{2-p}})$, and for some constant $K_4 = K_4(\theta, a_1)$.

\smallskip

Then let $\lambda_3 < \lambda_2$. One can prove that the first term in the right-hand side in \eqref{2termsp}  is bounded by
$\ds \exp \Big( - \frac{\lambda_3}{2} \, N \, \eps^2 \Big)$
provided
\begin{equation}\label{condRNeps}
R^2 \geq A \max ( 1, \eps^2, \ln (\eps^{-2} )) \, , \quad N \eps^2 \geq B \, 3^{C (R \eps^{-1})^{1/ \alpha}}
\end{equation}
for some positive constants $A, B$ and $C$ depending also on $\lambda_3$. Moreover, for
$\ds  \theta = \frac{\eps^{2-p} \, \lambda_3}{ 2 \, K_3} \raise 2pt \hbox{,}$
also the second term  in the right-hand side in \eqref{2termsp} is bounded by 
$\ds \exp \Big( - \frac{\lambda_3}{2} \, N \, \eps^2 \Big)$ as soon as
 $R^2 \geq R_4 \max ( 1, \ln (\eps^{-2} ))$, for some constant $R_4$ depending on $\lambda_3$.
 
\smallskip
 
Letting $\ds R = \eps \Big( \frac{1}{C \, \ln 3} \ln \frac{N \eps^2}{B} \Big)^{\alpha}$ if $\eps \in (0,1)$ and 
$R = \sqrt{A} \, \eps$ otherwise, and $\alpha' < \alpha$,  both conditions in \eqref{condRNeps} hold true as soon as 
$ N \geq N_0 \, \eps^{-2} \, \exp (N_0 \, \eps^{- 1 / \alpha'})$ for
some constant $N_0$ depending on $E_a, a, \lambda, \lambda_3, \alpha$ and $\alpha'$. 
Finally, given $\lambda' < \lambda_3 < \lambda$, this condition ensures that
\[
\pp\left[W_{p, [0, T]} (\mu,\hat{\mu}^N) > \eps \right] \leq 2 \,  \exp \Big( - \frac{\lambda_3}{2} \, N \, \eps^2 \Big) \leq \exp \Big( - \frac{\lambda'}{2} \, N \, \eps^2 \Big),
\]
possibly for some larger $N_0$. This concludes the argument in the case when $p \in [1,2)$.

\bigskip

In the case when $p=2$, given $0 < \eta < 1, \lambda_3 < \lambda_2$ and $a_2 < a_1$, the same condition on $N$ and $\eps$ (for some 
$N_0$) is sufficient for the bound
\[
\pp\left[W_{2, [0, T]} (\mu,\hat{\mu}^N) > \eps \right] \leq  \exp \Big( - \frac{\lambda_3}{2} \, \eta^2 \, N \, \eps^2 \Big) +
\exp \Big( - \frac{a_2}{2} (1 - \eta)^2 \, N \, \eps^2 \Big)
\] 
to hold (by \eqref{concl-trunc-2}). One optimizes this bound by letting
\[
a_2 = a\,  \frac{\lambda_3}{\lambda} \, (\in [0,a) ) \quad \text{and} \quad \eta = \frac{\sqrt{a_2}}{\sqrt{a_2} + \sqrt{\lambda_3}}
\cdot
\]
Given $\lambda' < \lambda_3 < \lambda$, this ensures the existence of $N_0$ such that
\[
\pp \big[W_{2, [0, T]} (\mu,\hat{\mu}^N) > \eps \big] \leq  2 \, 
\exp \Big( - \frac{\lambda_3}{2} \frac{a}{(\sqrt{a} + \sqrt{\lambda})^2} \, N \, \eps^2 \Big) \leq  
\exp \Big( - \frac{\lambda'}{2} \frac{1}{(1 + \sqrt{\lambda / a})^2} \, N \, \eps^2 \Big)
\]
for any $\eps >0$ and $N \geq N_0 \, \eps^{-2} \, \exp (N_0 \, \eps^{- 1 / \alpha'})$. This concludes the proof of Theorem
\ref{thmprelim} in this second and last case.

\medskip

\section{Coupling argument}\label{sectioncouplage}

Here begins the proof of Theorem \ref{thmprincipal} on the behaviour of our large interacting particle system.

We recall that we are given $N$ independent variables $X^i_{0}$ in $\rr^d$, with common law $\mu_0$, and $N$ independent Brownian motions
$B^i = (B^i_t)_{0 \leq t \leq T}$ in $\rr^d$, and we consider the solutions $X^i = (X^i_t)_{0 \leq t \leq T}$ to the coupled stochastic differential equations
\[
dX^i_t = \sqrt{2} \, dB^i_t - \nabla V(X^i_t) \, dt - \frac{1}{N} \sum_{j=1}^{N} \nabla W(X^i_t - X^j_t) \, dt, \qquad 1 \leq i \leq N . 
\]
We also let $\mu_{[0,T]}$ be the law of the process $Y = (Y_t)_{0 \leq t \leq T}$ defined by
\[
dY_t = \sqrt{2} \, dB_t - \nabla V (Y_t) \, dt - \nabla W * \mu_t (Y_t)\, dt 
\]
and starting at some $Y_0$ drawn according to $\mu_0$; here $B = (B_t)_{0 \leq t \leq T}$ also is a Brownian motion and $\mu_t$ is the
law of $Y_t$, that is, the time-marginal of $\mu_{[0,T]}$ at time $t$.

We want to compare this law $\mu_{[0,T]}$ and the empirical measure of the paths
\[
\hat{\mu}^N_{[0,T]} = \frac{1}{N} \sum_{i=1}^N \delta_{X^i}.
\]

For this purpose we introduce  $N$ independent processes $Y^i = (Y^i_t)_{0 \leq t \leq T}$ defined by
\begin{equation}\label{edsyicoup}
dY^i_t = \sqrt{2} \, dB^i_t - \nabla V (Y^i_t) \, dt - \nabla W * \mu_t (Y^i_t)\, dt, \qquad 1 \leq i \leq N
\end{equation}
for the same Brownian motions $B^i$, and such that $Y^i_0 = X^i_0$ initially. We let 
\begin{equation}\label{defnun}
\ds \hat{\nu}^N_{[0,T]} := \frac{1}{N} \sum_{i=1}^N \delta_{Y^i}
\end{equation}
and in this section we reduce the issue to measuring the distance between $\mu_{[0,T]}$ and $\hat{\nu}^N_{[0,T]}$.

\medskip

\begin{proposition}\label{propcouplage}
In the above notation and under the assumptions 
\[
\beta \, I \leq D^2 V (x) \, , \quad \gamma \, I \leq D^2 W(x) \leq \gamma' \, I \, , \qquad  x \in \rr^d
\]
on $V$ and $W$, where $\beta, \gamma$ and $\gamma'$ are real numbers, for any $T \geq 0$ there exists some constant $C$ depending only on $\beta, \gamma, \gamma'$ and $T$ 
such that 
\[
W_{1,[0, T]} (\mu_{[0,T]} , \hat{\mu}^{N}_{[0,T]}) \leq C \, W_{1,[0, T]} (\mu_{[0,T]} , \hat{\nu}^{N}_{[0,T]})
\]
almost surely. 
\end{proposition}

\smallskip

\begin{proof}
We first follow the lines of the proof of \cite[Proposition 5.1]{BGV05}, but in the end we want an estimate on 
the {\it trajectories}.
Since for each $i$ both processes $X^i$ and $Y^i$ are driven by the same Brownian motion $B^i$, the process $X^i - Y^i$ satisfies \[
d(X_t^i - Y_t^i) = - \bigl(\nabla V(X_t^i) - \nabla V(Y_t^i)\bigr) \, dt - 
\bigl(\nabla W \ast \hat{\mu}_t^N (X_t^i) - \nabla W\ast \mu_t(Y_t^i)\bigr) \, dt.
\]

In particular, letting $u \cdot v$ denote the scalar product of two vectors $u$ and $v$ in $\rr^d$, 
\begin{multline}\label{couplagecarre}
\frac{1}{2} \frac{d}{dt} \vert X_t^i - Y_t^i \vert^2 =
- \, \big(\nabla V(X_t^i) - \nabla V(Y_t^i) \big) \cdot (X_t^i - Y_t^i) \\ 
- \, \bigl(\nabla W \ast \hat{\mu}_t^N (X_t^i) - 
\nabla W\ast \mu_t(Y_t^i)\bigr) \, \cdot\,  (X_t^i - Y_t^i). 
\end{multline}
We  decompose the last term according to
$$
\nabla W \ast \hat{\mu}_t^N (X_t^i) - \nabla W\ast \mu_t(Y_t^i)= 
\big( \nabla W\ast\hat{\mu}_t^N-\nabla W\ast\mu_t \big) (X_t^i)
+ \big(\nabla W \ast \mu_t (X_t^i) - \nabla W\ast \mu_t(Y_t^i)\big).
$$
By our assumption on $D^2W$, the map $\nabla W(X_t^i - \, \cdot \,)$ is $\Gamma$-Lipschitz with $\Gamma := \max ( \vert \gamma \vert,
\vert \gamma' \vert)$. Consequently, by the Kantorovich-Rubinstein dual formulation \eqref{KR} of $W_{1, \tau}$,
$$
\Bigl| \nabla W \ast (\hat{\mu}_t^N -\mu_t)(X_t^i) \Bigr| = 
\left\vert \int_{\rr^d} \nabla W(X_t^i - y) \, d(\hat{\mu}_t^N - \mu_t)(y) \right\vert \leq 
\Gamma \, W_{1, \tau}(\hat{\mu}_t^N,\mu_t).
$$
Then \eqref{couplagecarre} and our convexity assumptions on $V$ and $W$ imply
$$
\frac{1}{2} \frac{d}{dt} \vert X_t^i - Y_t^i \vert^2  
\leq - (\beta + \gamma) \, \vert X_t^i - Y_t^i \vert^2
+ \Gamma \, W_{1, \tau}(\hat{\mu}_t^N,\mu_t) \, \vert X_t^i - Y_t^i \vert.
$$
In particular, by Gronwall's lemma,
$$
\vert X_t^i - Y_t^i \vert \leq \Gamma \int_{0}^{t} e^{-(\beta+\gamma)(t-u)} \, 
W_{1, \tau}(\hat{\mu}_u^N,\mu_u) \, du
$$
since initially $X^i_0 = Y^i_0$. Consequently, by convexity of the $W_{1, [0,t]}$ distance,
\begin{eqnarray}
W_{1, [0, t]}(\hat{\mu}^{N}_{[0,t]}, \hat{\nu}^{N}_{[0,t]}) 
& \leq & \frac{1}{N} \sum_{i=1}^{N} \sup_{0 \leq s \leq t} \vert X^i_s - Y^i_s \vert \nonumber\\
& \leq & \frac{1}{N} \sum_{i=1}^{N} \sup_{0 \leq s \leq t} \Gamma \int_{0}^{s} e^{-(\beta + \gamma)(s-u)}
W_{1, \tau}(\hat{\mu}^{N}_{u}, \mu_u) \, du \label{couplagegronwall}\\
& \leq & \Gamma e^{\vert \beta + \gamma \vert T} \int_{0}^{t} W_{1, \tau}(\hat{\mu}^{N}_{u}, \mu_u) \, du \nonumber
\end{eqnarray}
for all $0 \leq t \leq T$. But 
\[
W_{1, \tau}(\hat{\mu}^{N}_{u}, \mu_u) \leq 
W_{1, [0, u]}(\hat{\mu}^{N}_{[0,u]}, \mu_{[0,u]}) \leq 
W_{1, [0, u]}(\hat{\mu}^{N}_{[0,u]}, \hat{\nu}^{N}_{[0,u]})  + W_{1, [0, u]}(\hat{\nu}^{N}_{[0,u]}, \mu_{[0,u]}) 
\]
by the projection relation \eqref{ineq-projwp} and triangular inequality for $W_{1, [0,u]}$, so 
\begin{equation}\label{gronW}
W_{1, [0,t]}(\hat{\mu}^{N}_{[0,t]}, \hat{\nu}^{N}_{[0,t]}) 
\leq 
\Gamma e^{\vert \beta + \gamma \vert T} \int_{0}^{t} \exp \Big( \Gamma e^{\vert \beta + \gamma \vert T} (t-u) \Big) \, 
W_{1, [0, u]}(\hat{\nu}^{N}_{[0,u]}, \mu_{[0,u]}) \, du
\end{equation}
for all $0 \leq t \leq T$ by Gronwall's lemma again. Then, given $0 \leq u \leq t$,
\[
W_{1, [0, u]}(\hat{\nu}^{N}_{[0,u]}, \mu_{[0,u]}) \leq W_{1, [0, t]}(\hat{\nu}^{N}_{[0,t]}, \mu_{[0,t]})
\]
since $\hat{\nu}^{N}_{[0,u]}$ and $\mu_{[0,u]}$ are the respective image measures of 
$\hat{\nu}^{N}_{[0,t]}$ and $\mu_{[0,t]}$ by the $1$-Lipschitz map defined from $\CC([0,t], \rr^d)$ into $\CC([0,u], \rr^d)$ as the
restriction to $[0,u]$. Hence
\[
W_{1, [0, t]}(\hat{\mu}^{N}_{[0,t]}, \hat{\nu}^{N}_{[0,t]}) \leq C \, W_{1, [0, t]}(\hat{\nu}^{N}_{[0,t]}, \mu_{[0,t]})
\]
by \eqref{gronW}, for some constant $C$ depending only on $T, \beta, \gamma$ and $\gamma'$. This concludes the argument by triangular inequality.
\end{proof}

\medskip

\begin{remark}
If moreover $\beta + \gamma > \Gamma$, where again $\Gamma := \max ( \vert \gamma \vert, \vert \gamma' \vert)$ then we can let $C$ be $(\beta + \gamma) \, (\beta + \gamma - \Gamma)^{-1}$ in Proposition \ref{propcouplage}, independently of $T$. 
Indeed, if $\beta + \gamma >0$, then \eqref{couplagegronwall} leads to
\[
W_{1, [0, t]}(\hat{\mu}^{N}_{[0,t]}, \hat{\nu}^{N}_{[0,t]}) \leq \Gamma \sup_{0 \leq s \leq t} \int_{0}^{t} e^{-(\beta + \gamma)(s-u)}
\, du \; \sup_{0 \leq u \leq t} W_{1, \tau} (\hat{\mu}^{N}_{u}, \mu_u) \leq \frac{\Gamma}{\beta + \gamma} \, W_{1, [0, t]}  (\hat{\mu}^{N}_{[0,t]}, \mu_{[0,t]})
\]
and by triangular inequality
\[
W_{1,[0, t]}(\hat{\mu}^{N}_{[0,t]}, \mu_{[0,t]}) \leq \frac{\beta + \gamma}{\beta + \gamma - \Gamma} 
\, W_{1, [0, t]}(\hat{\nu}^{N}_{[0,t]}, \mu_{[0,t]})
\]
provided $\beta + \gamma > \Gamma$.

This is reminiscent of the fact that, under some convexity assumptions on $V$ and $W$, such as $\beta >0, \, \beta + 2 \, \gamma >0$, 
it has been proven in \cite{CMV03, CMV04, Mal01} that the time-marginal $\mu_t$ of the measure $\mu_{[0,t]}$ converges, as $t$ goes to
infinity, to the stationary solution to the limit equation \eqref{edp}. One can also prove in this context that (in expectation) 
observables of the particle system are bounded in time. 

Hence, under this kind of assumptions, one could hope for some uniform in time constants in this coupling argument: that was obtained in
\cite[Proposition 5.1]{BGV05} for the time-marginals, and here for the whole processes. However, contrary to \cite{BGV05} where this property was used to approach the stationary solution by coupling together estimates of concentration of the empirical measure (as $N$ goes to infinity) with estimates of convergence to equilibrium (as $t$ goes to infinity), in this work we are concerned
with finite time intervals only, and shall not use this specific property in the sequel.
\end{remark}

\medskip

\section{Integrability in H\"older norm}\label{sectionconclusion}

In the previous section we have reduced the issue of measuring the distance between  $\mu_{[0,T]}$ and $\hat{\mu}^N_{[0,T]}$ to measuring
the distance between $\mu_{[0,T]}$ and the empirical measure $\hat{\nu}^{N}_{[0,T]}$ of $N$ independent random variables drawn according 
to $\mu_{[0,T]}$.

We now solve the latter issue by proving that the measure $\mu_{[0,T]}$ fulfills the hypotheses of Theorem \ref{thmprelim} with
$p = 1$, namely, that there exist $\alpha \in (0,1]$ and $a > 0$ such that
\[
\int_{\CC} e^{a \Vert x \Vert_{\alpha}^2} \, d\mu_{[0,T]} (x) := \ee \exp ( a \Vert Y \Vert_{\alpha}^2 ) < + \infty.
\]
Here again $\CC$ stands for ${\mathcal C}([0,T], \rr^d)$, $\Vert f \Vert_{\alpha}$ for the H\"older norm
of a function $f$ on $[0,T]$, and $Y = (Y_t)_{0 \leq t \leq T}$ is the solution to the stochastic differential equation
\begin{equation}\label{sdeconclusion}
dY_t = \sqrt{2} \, dB_t - \nabla V (Y_t) \, dt - \nabla W * \mu_t (Y_t) \, dt
\end{equation}
starting at $Y_0$ drawn according to $\mu_0$, where $\mu_t$ is the law of $Y_t$. 

\begin{proposition}\label{propholder}
Let $\mu_0$ be a probability measure on $\rr^d$ admitting a finite square-exponential moment and let $Y_0$ be drawn according to $\mu_0$. Given $T \geq 0$, $V$ and $W$ satisfying
hypotheses \eqref{D2V,D2W}, let $Y$ be the solution to \eqref{sdeconclusion} starting at $Y_0$. Then, for any $\alpha \in (0, 1/2)$, there exists $a > 0$, depending on $\mu_0$ only through a finite square-exponential moment, such that $\ee \exp ( a \Vert Y \Vert_{\alpha}^2) $ be finite.
\end{proposition}

\smallskip

Assuming this result for the moment we can now conclude the {\it proof of Theorem \ref{thmprincipal}}. Let indeed $\alpha$  
be given in $(0, 1/2)$, and $\alpha_0 \in (\alpha, 1/2)$. Then, by Proposition \ref{propholder} and Theorem \ref{thmprelim}, applied
with $\alpha = \alpha_0$ and $\alpha' = \alpha$, there exist some constants $\tilde{K}$ and $\tilde{N}_0$, depending on $\alpha_0, \alpha, T$ and a square-exponential 
moment of $\mu_0$, such that
\[
\pp \, \big[ W_{1, [0,T]} ( \mu_{[0,T]} , \hat{\nu}^{N}_{[0,T]} ) > \tilde{\eps} \big] \leq e^{ - \tilde{K} N \tilde{\eps}^2}
\]
for any $\tilde{\eps} > 0$ and $N \geq \tilde{N}_0 \, \tilde{\eps}^{-2} \, \exp(\tilde{N}_0 \, \tilde{\eps}^{- 1/ \alpha})$, where $\hat{\nu}^{N}_{[0,T]}$ is defined by \eqref{edsyicoup} and \eqref{defnun}. Then, 
by Proposition \ref{propcouplage}, there exist some constants $C$, depending only on $T$, and then $K$ and $N_0$, 
depending on $\alpha_0, \alpha, T$ and a finite square-exponential moment of $\mu_0$, such that
\[
\pp \, \big[ W_{1,[0,T]} ( \mu_{[0,T]} , \hat{\mu}^{N}_{[0,T]} ) > \eps ] \leq \pp \, \big[ W_{1,[0,T]} ( \mu_{[0,T]} , \hat{\nu}^N_{[0,T]} ) 
> \eps / C ]  \leq e^{ - K N \eps^2 }
\]
for any $\eps > 0$ and $N \geq N_0 \, \eps^{-2} \, \exp( N_0 \, \eps^{- 1 / \alpha})$.
This concludes the argument. $\hfill \Box$

\medskip

\begin{proof}[Proof of Proposition \ref{propholder}]
It is necessary and sufficient to prove that there exist positive cons\-tants $a_1$ and $a_2$ such that
$\ee \exp ( a_1 \Vert Y \Vert_{\infty}^2 )$ and $\ee \exp ( a_2 [Y]_{\alpha}^2 )$ be finite, where 
$[ \, \cdot \, ]_{\alpha}$ stands for the H\"older seminorm defined in Section \ref{subsectionprelim}.

\smallskip

{\bf Step 1.}  We start with the expectation in uniform norm. For this we first note that, according to
\cite[Proposition 3.1]{BGV05}, there exist some positive constants $M$ and $\overline{a}$, depen\-ding on $\mu_0$ only through a finite square-exponential moment, such that
$\ds \sup_{0 \leq t \leq T} \ee \, \vert Y_t \vert^2$  and
$\ds \sup_{0 \leq t \leq T} \ee \exp ( \overline{a} \vert Y_t \vert^2 )$ be finite and bounded by $M$.

Then we let $b$ be some smooth function on $[0,T]$, to be chosen later on, 
and we let $Z_t = \exp  ( b(t) \, \vert Y_t \vert^2)$. We want to prove that $\ds \ee \, \sup_{0 \leq t \leq T}  Z_t$ is finite for some positive 
function $b$. By It\^o's formula,
\[
Z_t = Z_0 + M_t + \int_{0}^{t} \Big[ b'(s) \vert Y_s \vert^2 + 2 \, d \, b(s) + 4 \, b(s)^2 \, \vert Y_s \vert^2
- 2 \, b(s) \, Y_s \, \cdot \, (\nabla V (Y_s) + \nabla W * \mu_s (Y_s)) \Big] \, Z_s  \, ds
\]
where $(M_t)_{0 \leq t \leq T}$ is the martingale defined as
\[
M_t = 2 \sqrt{2} \int_{0}^{t} b(s) \, Z_s \, Y_s \, \cdot \, dB_s.
\]
But $D^2 V(x) \geq \beta \, I$ for all $x \in \rr^d$, so for any $\delta >0$ and $y \in \rr^d$ we have
\[
- y \, \cdot \, \nabla V (y)  \leq 
(\delta - \beta ) \vert y \vert^2 + \frac{\vert \nabla V(0) \vert^2}{4 \delta} \cdot
\]
Furthermore $\nabla W$ is $\Gamma$-Lipschitz and $\nabla W(0) = 0$, so
$$
 - 2 \,  y \, \cdot \, \nabla W * \mu_s (y) 
 \leq 2 \, \Gamma
\int_{\rr^d} \vert y \vert \, \vert y-z \vert \, d\mu_s(z) \leq 3 \, \Gamma \, \vert y \vert^2 + \Gamma \int_{\rr^d} \vert z \vert^2 \, d\mu_s(z).
$$
But $\ds \int_{\rr^d} \vert z \vert^2 \, d\mu_s(z) = \ee \, \vert Y_s \vert^2 \leq M$ on $[0,T]$, so collecting
all terms together, we obtain
\[
Z_t \leq Z_0 + M_t + \int_{0}^{t} \big[ C(s) + D(s) \, \vert Y_s \vert^2 \big] \, Z_s \, ds
\]
where $ C(s) =  \ds \Big(  2 \, d + \Gamma \, M +  \frac{\vert \nabla V(0) \vert^2}{2 \delta} \Big) \, b(s)$ and 
$D(s) = b'(s) + 4 \, b(s)^2 + \big( 2 \, (\delta - \beta) + 3 \, \Gamma \big) \, b(s).$
Given $\delta >0$ such that $c := 2 \, (\delta - \beta) + 3 \, \Gamma$ be positive, we let $b(s)$ such that $D(s) \equiv 0$,~namely 
$$
b(s) = e^{-c s} \left( b(0)^{-1} + 4\,  c^{-1} (1 - e^{- c s}) \right)^{-1}
$$
for some $b(0)$ to be chosen later on. In particular $b$ is a nonincreasing continuous positive function on $[0, +\infty)$, and,
 for this function $b$, $Z_t$ almost surely satisfies the inequality
\[
Z_t \leq Z_0 + M_t + C(0) \int_{0}^{t} Z_s \, ds .
\]
In particular
\begin{equation}\label{espsupzt}
\ee \sup_{0 \leq t \leq T} Z_t \leq \ee \, Z_0 + \ee \sup_{0 \leq t \leq T} M_t + C(0) \int_{0}^{T} \ee \, Z_s \, ds .
\end{equation}

But, by Cauchy-Schwarz' and Doob's inequalities,
\[
\Big( \ee \sup_{0 \leq t \leq T} M_t \Big)^2 \leq \ee \sup_{0 \leq t \leq T}  \vert M_t \vert^2 
\leq 2 \,  \sup_{0 \leq t \leq T} \ee \, \vert M_t \vert^2.
\]
Then, by It\^o's formula again,
\begin{eqnarray*}
\ee \vert M_t \vert^2
 =  8 \int_{0}^{t} b(s)^2 \, \ee\,  \big[ Z_s^2 \, \vert Y_s \vert^2 \big] \, ds &\leq & 8  \, b(0) \int_{0}^{t} \ee \, \Big[ b(s) \, \vert Y_s \vert^2 \, \exp( 2 \, b(s) \, \vert Y_s \vert^2 ) \Big] \, ds \\
& \leq & 8 \, b(0)\int_{0}^{t} \ee \, \exp( 3 \, b(0) \, \vert Y_s \vert^2 ) \, ds.
\end{eqnarray*}
Choosing $b(0) \leq \overline{a} / 3$, this ensures that
$ \ds \sup_{0 \leq t \leq T} \ee \vert M_t \vert^2 $, whence $\ds \ee \sup_{0 \leq t \leq T} M_t$, is finite.

\smallskip

Since, for this $b(0)$, $\ds \sup_{0 \leq t \leq T} \ee \, Z_t$ also is finite, it follows from \eqref{espsupzt} that
so is $\ds \ee \sup_{0 \leq t \leq T} Z_t$, which concludes the argument for the expectation in uniform norm with $a_1 = b(T)$. 

\bigskip

{\bf Step 2.} We now turn to the expectation in H\"older seminorm. Writing the solution as
\[
Y_t = Y_0 + B_t - \int_{0}^{t} \big( \nabla V (Y_s) + \nabla W * \mu_s (Y_s) \big) \, ds
\]
we obtain
\[
[Y ]_{\alpha} \leq [B]_{\alpha} + 
\Bigl[ \int_{0}^{{\text{\bf{.}}}} \big( \nabla V (Y_s) + \nabla W * \mu_s (Y_s) \big) \, ds \Big]_{\alpha}
\]
almost surely; here $Y$ and $B$ stand as before for the map $t \mapsto Y_t$ and $t \mapsto B_t$ 
respectively, and 
$\ds \int_{0}^{{\text{\bf{.}}}} \varphi(s) \, ds$ is an antiderivative of $\varphi$. Hence, by Cauchy-Schwarz' inequality,
\[
\ee \exp (a_2 [Y]_{\alpha}^2) \leq \big( \ee \exp (4 \, a_2 [B]_{\alpha}^2) \big)^{1/2}
\Big( \ee \exp \Big( 4 \, a_2 \Bigl[ \int_{0}^{{\text{\bf{.}}}} ( \nabla V (Y_s) + \nabla W * \mu_s (Y_s) ) \, ds \Bigr]_{\alpha}^2 \Big) \Big)^{1/2}.
\]

But, on one hand, $\ds \ee \, \exp \big( 4 \, a_2 [ B ]_{\alpha}^2 \big)$ is finite for $a_2$ small enough (see 
\cite[Theorem 1.3.2]{Fer75} for instance, with $E = \CC$ and $N(f) = [ f ]_{\alpha})$. On the other hand, by assumption \eqref{D2V,D2W},
$\nabla V$ and $\nabla W$ are respectively B- and $\Gamma$-Lipschitz with B $:= \max (\vert \beta \vert, \vert \beta' \vert)$ and $\Gamma := \max (\vert \gamma \vert, \vert \gamma' \vert)$, so there exists some constant A such that
\[
\big\vert \nabla V (y) + \nabla W * \mu_s (y) \big\vert \leq \text{A} + (\text{B} + \Gamma) \vert y \vert\, 
\]
for all $y \in \rr^d$ and $s \in [0, T]$. In particular
\begin{eqnarray*}
\Bigl[ \int_{0}^{{\text{\bf{.}}}} \big( \nabla V (Y_s) + \nabla W * \mu_s (Y_s) \big) \, ds \Bigr]_{\alpha}
& \leq &  
\sup_{0 \leq s , t \leq T} \frac{1}{\vert t-s \vert^{\alpha}}  \int_{s}^{t} \big( \text{A} + (\text{B} + \Gamma) \vert Y_u \vert \big)\, du  \,
 \\
& \leq &   T^{1 -  \alpha} \big( \text{A} + (\text{B} + \Gamma) \Vert Y \Vert_{\infty} \big) 
\end{eqnarray*}
almost surely, and
\begin{multline*}
\ee \exp \left( 4 \, a_2 \Big[ \int_{0}^{{\text{\bf{.}}}} \big( \nabla V (Y_s) + \nabla W * \mu_s (Y_s) \big) \, ds \Bigr]_{\alpha} ^2 \right)\\
\leq \exp \big( 8 \, a_2 \, T^{2-2 \alpha} \text{A}^2 \big) \;  
\ee \exp \Big( 8 \, a_2 \, T^{2 - 2 \alpha} (\text{B} + \Gamma)^2 \Vert Y \Vert_{\infty}^2 \Big) 
\end{multline*}
which by step 1 is finite as soon as $8 \,  a_2 \, T^{2 - 2 \alpha} (\text{B} + \Gamma)^2 \leq a_1$. 

\smallskip

On the whole, $\ee \exp ( a_2 [Y]_{\alpha}^2)$ is indeed finite for $a_2$ small enough, depending on $\mu_0$ only through a finite square-exponential moment, which concludes the argument.
\end{proof}

\medskip

\appendix

\section*{Appendix. Metric entropy of a H\"older space}\label{secannexe}

 The aim of this appendix is
to establish the bound \eqref{cov3} used in the covering argument in the proof of Theorem \ref{thmprelim},  which amounts to studying the metric entropy of a H\"older space and of some related space of probability measures.

\bigskip

In the notation introduced in Sections \ref{subsectionnotation} and \ref{subsectionprelim}, it follows from Ascoli's theorem 
that the closed ball ${\mathcal B}^\alpha_R := {\mathcal B}^\alpha_R([0,T], \R^d) = \{ f \in \CCalpha ; \Vert f \Vert_{\alpha} \leq R \}$ of 
center $0$ and radius $R$ in $\CCalpha$  is a compact metric space for the metric defined by the uniform norm.
Here we estimate by how many balls of given radius $r < R$ and centered in ${\mathcal B}^\alpha_R$ 
the compact metric space  ${\mathcal B}^\alpha_R$
can be covered. We note that for  $r \geq R$  the sole ball $\{ f \in {\mathcal B}^\alpha_R; \Vert f \Vert_{\infty} \leq r \}$
covers   ${\mathcal B}^\alpha_R$.

\bigskip

{\bf Notation:} Given $r >0$, the {\em covering number} 
${\mathcal N}(E, r)$ of a compact metric space $(E,d)$ is defined as the infimum of the integers $n$ such that $E$ can be covered by $n$ balls centered in $E$ and of radius $r$ in $d$ metric. 
Then  we have  the following result which gives some lower and upper bounds on the covering number 
${\mathcal N}({\mathcal B}^\alpha_R,r)$ and in our case makes more precise the bounds given for instance in \cite{Lor66} or \cite{VW95}:

\medskip

\begin{theorem}\label{proposition2}
Given some integer number $d \geq 1$,  some positive  numbers $T$, $R$, $r$ and $\alpha $ with $r<R$ and $\alpha \leq 1$,
the covering number ${\mathcal N}({\mathcal B}^\alpha_R,r)$ of ${\mathcal B}^\alpha_R$, equipped
 with the uniform norm,  satisfies
\[
{\mathcal N}({\mathcal B}^\alpha_R,r) 
\leq 
\Big(10 \, \sqrt{d} \, \frac{R}{r} \Big)^d \,3^ {5^\frac{1}{\alpha}\,d^{1 + \frac{1}{2 \alpha}} \,T\,(\frac{R}{r})^\frac{1}{\alpha}}.
\]
If  moreover, for instance, $\displaystyle r  \leq \frac{T^\alpha}{4 \, T^\alpha+4}R$, then 
\[
{\mathcal N}({\mathcal B}^\alpha_R,r)\, \geq \, \Big( \frac{\sqrt{d}}{4}\,\frac{R}{r} \Big)^d \, 
2^{2^{-\frac{1}{\alpha}}\,d^{1 + \frac{1}{2 \alpha}} \, T\, (\frac{R}{r})^{\frac{1}{\alpha}}} .
\]
\end{theorem}

\medskip

The lower bound ensures that the upper bound, from which depends the condition on the size of the sample in Theorems \ref{thmprelim} and hence \ref{thmprincipal}, has the good order of growth in $R/r$.

\medskip
\begin{proof}  
{\bf 1.} We start by establishing the {\bf upper bound}. 
\smallskip

1.1. We first consider the case  when $d=1$. 
\smallskip

Given  $J$ and  $K$  some  integers larger or equal to $1$, we let $\displaystyle \tau=\frac{T}{J}$ and 
$\displaystyle \eta=\frac{R}{K}\raise 2pt \hbox{,}$ and then
$$
\begin{array}{rcllrl}
t_j & = & (j-\frac{1}{2}) \, \tau \, , \qquad & j \in \nn  \,,  \quad & 1  & \leq j \leq J, \\
y_k & = & (k-\frac{1}{2}) \, \eta \, , & k \in \nn \, , & -K+1  & \leq k \leq K.
\end{array}
$$
Then we cover the rectangle $[0,T]\times [-R,+R]$ in $\R_t \times \R_y$, which  contains the graph of all functions in 
${\mathcal B}^\alpha_R([0,T], \R)$, by a lattice with step  $\tau$ in $t-$axis and $\eta$ in  $y-$axis.

Then let  $f$ be a given  function in ${\mathcal B}^\alpha_R([0,T], \R) $. Since the  intervals 
$\displaystyle  [y_k-\frac{\eta}{2}\,\raise 2pt \hbox{,}\, y_k+\frac{\eta}{2}]$ cover the interval $[-R,+R]$, for every 
integer $j \in [1,J]$ there exists some integer $k(j) \in [-K+1,+K]$ such that
\[
\vert f(t_j)-y_{k(j)} \vert \leq \frac{\eta} {2}\cdot
\]
In particular
\[
\vert y_{k(j+1)}-y_{k(j)}\vert \leq \frac{\eta}{2}+ \vert f(t_{j+1})-f(t_{j})\vert +\frac{\eta}{2}
\leq \eta+R \, \vert t_{j+1}-t_{j}\vert^\alpha 
\leq \eta+ R \,\tau^\alpha < 2\eta
\]
if  we suppose $K T^\alpha<J^\alpha.$ But since the  $y_k$ take values which are regularly distant of  $\displaystyle \eta,$ 
it follows that more precisely
\[
\vert y_{k(j+1)}-y_{k(j)}\vert\leq \eta .
\]

From this map $\ k:[1,J] \cap \N \to [-K+1,K]\cap \N$, we define the function
$f_k:[0,T] \to [-R,+R]$ affine on each interval of the subdivision $(0,t_1,\cdots,t_J,T)$ and such that
\begin{eqnarray*}
f_{ k}(0) & = & f_{ k}(t_1), \\
f_k(t_j)  & = & y_{ k(j)}, \quad   1 \leq  j  \leq J \\
f_{k}(T)  & = & f_{ k}(t_J).
\end{eqnarray*}
In particular we note that  this function $f_k$ is Lipschitz with
\[
\sup_{0 \leq t,s \leq T}\frac {\vert f_k (t)-f_k(s)\vert }{\vert t-s \vert}
= \sup_{1 \leq k \leq K}\frac {\vert y_{k(j+1)}-y_{k(j)}\vert }{\vert t_{j+1}-t_j\vert}\,
\leq \,\frac{\eta}{\tau}
\]
but that it does not necessarily belong to ${\mathcal B}^\alpha_R([0,T], \R) $.

The number of such functions $f_k$  is bounded by the number of  $J$-uples $(y_{k(j)})_{1 \leq j \leq J}$ such that 
$\displaystyle \vert y_{k(j+1)}-y_{k(j)}\vert \leq \eta$ for  $1 \leq j \leq J-1$, that is the number of   $J$-uples 
$(k(j))_{1 \leq j \leq J}$ such that $\vert k(j+1)-k(j) \vert \leq 1$ for $1 \leq j \leq J-1$. Such $J$-uples are 
obtained by choosing  $k(1)$ among $2K$ values, then 
$k(2)$ among $3$ values for $-K+2\leq k(1)\leq +K-1$ or $2$ values for $k(1)=-K+1$ and $+K$, and so on. Hence  
there exist at most $ 2 \, K \, 3^{J-1}$ such functions $f_k$.

\smallskip

If we now let $K$  be the smallest integer larger or equal to $\displaystyle 4\,\frac{R}{r}$ and  $J$ such that  
$K T^\alpha<J^\alpha$, then
\[
\Vert f -f_k \Vert_\infty \leq \frac {r}{2} \cdot
\]
Indeed, given  $t$ in $  [0,T]$, there exists some integer number $j $ in $[1,J]$ such that $t$ belongs to 
$\displaystyle  [t_j-\frac{\tau}{2}\,\raise 2pt \hbox{,}\,t_j+\frac{ \tau}{2}]$, so that
\begin{multline*}
\vert f(t) -f_k(t) \vert \leq \vert f(t) -f(t_j)\vert +\vert f(t_j) -f_k(t_j) \vert +\vert f_k(t_j) -f_k(t)\vert
 \\
\leq R \, \vert t-t_j\vert^\alpha +\vert f(t_j) -y_{j(k)} \vert +\frac{\eta}{\tau}\, \vert t_j-t \vert 
\leq R \Big( \frac {\tau}{2} \Big)^\alpha +\frac {\eta}{2}+\frac{\eta}{\tau} \frac {\tau}{2} \leq 2\eta \leq \frac {r}{2}\cdot
\end{multline*}

Hence we can cover ${\mathcal B}^\alpha_R([0,T], \R) $ by less than $ 2 \, K \, 3^{J-1}$
balls of radius $\displaystyle \frac{r}{2}$ of the metric space  ${\mathcal C}([0,T], \R) $
equipped with the   uniform norm, and if we let  $J$ and  $K$ be the smallest integers larger or equal to   
$\displaystyle 5^\frac{1}{\alpha} \, T\,(\frac{R}{r})^{\frac{1}{\alpha}}$ and $\displaystyle 4\,\frac{R}{r}$  respectively, then
$K T^\alpha<J^\alpha$ holds true and
\[
2 \, K \, 3^{J-1} \leq 10 \,\frac{R}{r}\,3^ {5^\frac{1}{\alpha}\,T\,(\frac{R}{r})^\frac{1}{\alpha}}.
\]
\medskip

1.2.  From this we now deduce the upper bound in the general case $d \geq 1$.

\smallskip

Let $F $ be a given  function in  ${\mathcal B}^\alpha_R([0,T], \R^d) $ with components 
$F_i \in {\mathcal B}^\alpha_R([0,T], \R) $ for $1 \leq i \leq d$. Let now $J$ and  $K$ be the smallest integers larger or equal to   
$\displaystyle 5^\frac{1}{\alpha} \, T\,(\sqrt{d} \, \frac{R}{r})^{\frac{1}{\alpha}}$ and $\displaystyle 4\, \sqrt{d} \, \frac{R}{r}$ respectively. With each $i$, we associate an  integer $k_i$ in  
$[1, 2 \, K \, 3^{J-1}]$ such that 
\[
\Vert F_i-f_{k_i} \Vert_\infty \leq \frac{ r}{2 \, \sqrt{d}}
\]
where the $f_k$ are the functions in $\CC([0,T], \rr)$ defined in the first step (relatively to $\ds \frac{r}{\sqrt{d}}$ instead of $r$).

Then the function  $F_{k_1,\cdots ,k_d}$ with components $f_{k_i}$
for $1 \leq i \leq d$  belongs to  $ {\mathcal C}([0,T], \R^d)$~and satisfies 
$\displaystyle \Vert F-F_{k_1,\cdots, k_d} \Vert_\infty \leq \frac{r}{2} \cdot$ Moreover there are at most $(2 \, K \, 3^{J-1} )^d$
such functions~$F_{k_1, \cdots, k_d}$.

Consequently we can cover  ${\mathcal B}^\alpha_R([0,R],\R^d) $ by less than $(2 \, K \, 3^{J-1})^d$ balls of radius  
$\displaystyle \frac {r}{2}$ of the metric space $ {\mathcal C}([0,T], \R^d)$ equipped with  the uniform norm,
whence by less than $\big(2 \, K \, 3^{J-1} \big)^d$
balls of radius  $r$ of the metric space   ${\mathcal B}^\alpha_R([0,T], \R^d)$ equipped with  the uniform norm. 

\medskip

This concludes the proof of the upper bound of the covering number   ${\mathcal N}({\mathcal B}^\alpha_R([0,T], \R^d),r)$.

\bigskip

\smallskip

{\bf 2.} We now turn to the {\bf lower bound}.

\smallskip

2.1. We first consider the case $d=1$.

We can give different types of lower bounds by considering special functions of the type $f_k$ defined in the first step.  Here, for instance, we give the detail for one of them.
 
\smallskip

Given some non-zero integer $J$, we let $\displaystyle \tau=\frac{T}{J}$ and 
$\displaystyle \eta=\tau^\alpha \, R$, and then
$$
\begin{array}{rcllrl}
\ds t_j & = & (j-\frac{1}{2}) \, \tau \, , \qquad  & j \in \nn \, , \qquad  & 1 & \leq j \leq J, \\
\ds y_k & = & (k-\frac{1}{2}) \, \eta  \, , & k \in \nn \, , & -\tau^{-\alpha}+\frac{1}{2}  & \leq  k \leq \tau^{-\alpha}+\frac{1}{2}\cdot
\end{array}
$$

From a map $\ k:[1,J] \cap \N \to [0,1]\cap \N$, we define as above the function
$f_k:[0,T] \to [y_0,y_1]$ affine on every interval of the subdivision $(0,t_1,\cdots,t_J,T)$ and such that
\begin{eqnarray*}
f_{ k}(0) & = & f_{ k}(t_1) \\
f_k(t_j)  & = & y_{ k(j)}, \qquad   1 \leq  j  \leq J \\
f_{k}(T)  & = & f_{ k}(t_J).
\end{eqnarray*}

Given some integer $\ell$ such that $\displaystyle-\tau^{-\alpha}+\frac{1}{2} \leq \ell \leq \tau^{-\alpha}-\frac{1}{2} \raise 2pt \hbox{,}$ we 
define the function $f_{k \ell}:[0,T] \to [y_{\ell},y_{\ell+1}]$ such that
\[
f_{k \ell}(t)=f_{k}(t)+ \ell \, \eta.
\]
Then  $f_{k \ell}$ belongs to ${\mathcal B}^\alpha_R([0,T], \R)$ and 
$\displaystyle \Vert f_{k \ell} -  f_{k' \ell'} \Vert_\infty \geq \eta$ if $f_{k \ell} \neq f_{k' \ell'}$.

If for instance $r  <  \inf(R,2^{-1}T^\alpha R)$ and $J+1$ is the smallest integer larger 
or equal to $\displaystyle2^{-\frac{1}{\alpha}}T\,(\frac{R}{r})^{\frac{1}{\alpha}}$, then 
$\displaystyle \Vert f_{k \ell} -  f_{k' \ell'} \Vert_\infty > 2r$ if $ f_{k \ell} \neq f_{k' \ell'}$.

Thus we have found $L\,2^J$ elements in ${\mathcal B}^\alpha_R([0,T], \R)$ mutually distant of at least
$2r$ in uniform norm, where $L$ is the number of integers $ \ell$ between 
 $\displaystyle -\tau^{-\alpha}+\frac{1}{2}$ and $\ds  \tau^{-\alpha}-\frac{1}{2}\cdot$ Hence
\[
{\mathcal N}({\mathcal B}^\alpha_R([0,T], \R),r) \,\geq \,L\,2^J .
\]
But
\[
L\,>\, 2 \big( (\tau^{-\alpha}-\frac{1}{2})-1\big)+1 \,= \, 2\tau^{-\alpha}-2\,\geq \,
\big( (\frac{R}{r})^{\frac{1}{\alpha}}-\frac{2^{\frac{1}{\alpha}}}{T} \big)^\alpha-2\, \geq \,\frac{R}{r}-\frac{2}{T^\alpha} -2.
\]
If moreover, for instance, $\displaystyle r  \leq \frac{T^\alpha}{4 \, T^\alpha+4}R$, then $\displaystyle L\geq \frac{R}{2r}$ and 
 \[
{\mathcal N}({\mathcal B}^\alpha_R([0,T], \R),r)\, \geq \, \frac{1}{4}\, \frac{R}{r} \,
2^{2^{-\frac{1}{\alpha}}\,T\, (\frac{R}{r})^{\frac{1}{\alpha}}}.
\]

\medskip

2.2. From this we now deduce the lower bound in the general case $d \geq 1$.
\smallskip

The $L^d\,2^{dJ} $ functions $F_{k_1 \ell_1,\cdots ,k_d \ell_d}$ with components $f_{k_j \ell_j}$ for $j=1,\cdots,d$ where  $f_{k_j \ell_j}$  
have been defined in the step 1, belong to ${\mathcal B}^\alpha_R([0,T], \R^d)$ and are mutually distant of at least  $2 \, \sqrt{d} \, r$. 

This concludes the argument for the lower bound of the  number   ${\mathcal N}({\mathcal B}^\alpha_R([0,T], \R^d),r)$.
\end{proof}

\bigskip

We now turn to the covering number of the corresponding space of probability measures: given a Polish metric 
space $(E, d)$, $p \geq 1$ and $\delta >0$, we denote ${\mathcal N}_{p} ({\mathcal P}(E), \delta)$ the 
covering number of ${\mathcal P}(E)$ for the $W_p$ distance.

\bigskip

Then we have the following general result which is proven in \cite{BGV05} (see also \cite{DZ98}, \cite{KZ95}):

\begin{theorem}
Let $(E,d)$ be a Polish metric space with finite diameter $D$, and let $p$ and  $\delta $ be some real numbers  with $p \geq 1$ and  
$0 < \delta < D$.  Then the covering number ${\mathcal N}_p({\mathcal P}(E),\delta)$ of ${\mathcal P}(E)$ satisfies
\[
{\mathcal N}_p({\mathcal P}(E),\delta) \leq \big(8\,e\, \frac{D}{\delta}\big)^ {p\,{\mathcal N}(E,\frac{\delta}{2})}.
\]
\end{theorem}

\bigskip
Note that if $\delta \geq D$ we simply have ${\mathcal N}_p({\mathcal P}(E),\delta) =1$ since  the Wasserstein distance  
between any two  probability measures  on $E$ is at most $D$.
\bigskip

Since  ${\mathcal B}^\alpha_R$  equipped with the  metric defined by the uniform norm is a Polish metric space  
with finite diameter $2R$,  we deduce  the following result:
\bigskip

\begin{theorem}\label{thm-annexe}
Let $d \geq 1$, $p$, $T$, $R$,  $\delta$ and  $\alpha $ be some positive numbers with $p \geq 1$,  
$\delta <2R$ and  $\alpha  \leq 1$. Let also ${\mathcal B}^\alpha_R = \{f \in \CCalpha ; \Vert f \Vert_{\alpha} \leq R \}$ be equipped with  the uniform norm. Then  
the space ${\mathcal P}({\mathcal B}^\alpha_R)$ of probability measures on ${\mathcal B}^\alpha_R$ 
can be covered by ${\mathcal N}_p({\mathcal P}({\mathcal B}^\alpha_R),\delta)$ balls of radius $\delta$ in Wasserstein distance $W_p$, with
\[
{\mathcal N}_p({\mathcal P}({\mathcal B}^\alpha_R),\delta) 
\leq 
\big(16\,e\, R \delta^{-1}\big)^{p\, (20\, \sqrt{d} \,R\,  \delta^{-1})^d\,
3^ {{10}^\frac{1}{\alpha}d^{1 + \frac{1}{2 \alpha}} \, T\, (R\delta^{-1})^\frac{1}{\alpha}}}.
\]
For $\delta \geq 2R$, we have 
$$
{\mathcal N}_p({\mathcal P}({\mathcal B}^\alpha_R),\delta)=1.
$$
\end{theorem}

 \bigskip

\noindent {\bf Acknowledgments.} The author thanks Professors M. Ledoux for his kind interest in this work, in particular pointing out Reference \cite{Fer75}, and A. Guillin and C. Villani for stimulating discussions 
during the preparation of \cite{BGV05}.

\end{document}